\documentclass[graybox]{svmult}

\usepackage{mathptmx}       
\usepackage{helvet}         
\usepackage{courier}        
\usepackage{makeidx}         
\usepackage{graphicx}        
\usepackage{multicol}        
\usepackage[bottom]{footmisc}

\usepackage{cite}
\usepackage{fancyvrb}        
\usepackage{times}
\usepackage{amssymb}
\usepackage{amsmath}
\usepackage{caption}

\newcommand{\Div}{\mathop{\operatorname{div}}\nolimits}
\newcommand{\VE}{{\mathbf{E}}}
\newcommand{\VD}{{\mathbf{D}}}
\newcommand{\Vx}{{\mathbf{x}}}
\newcommand{\Vy}{{\mathbf{y}}}
\providecommand{\Vn}{{\mathbf{n}}}
\providecommand{\grad}{\operatorname{{\bf grad}}}

\begin{document}

\title*{Dielectric Breakdown Prediction with GPU-Accelerated BEM}
\author{Cedric M\"unger 
\and
Steffen B\"orm 
\and
J\"org Ostrowski 
}
\institute{Cedric M\"unger \at Seminar for Applied Mathematics, ETH Z\"urich, Switzerland,
\email{cmuenger@student.ethz.ch}
  \and Steffen B\"orm \at Math. Seminar, Christian-Albrechts-Univ. Kiel, Germany,
  \email{boerm@math.uni-kiel.de}
  \and J\"org Ostrowski \at ABB Corporate Research, Baden, Switzerland,\email{joerg.ostrowski@ch.abb.com}}

%
%
\maketitle

\abstract{The prediction of a dielectric breakdown in a high-voltage device is based on criteria that evaluate the electric field along possible breakdown paths. For this purpose it is necessary to efficiently compute the electric field at arbitrary points in space. A boundary element method (BEM) based on an indirect formulation, realized with MPI-parallel collocation, has proven to cope very well with this requirement. It deploys surface meshes only, which are easy to generate even for complex industrial geometries. The assembly of the large dense BEM-matrix, as well as the iterative solution of the resulting system, and the evaluation along the field-lines all require to carry out the same type of calculation many times. Graphical Processing Units (GPUs) promise to be more efficient than Central Processing Units (CPUs) when it is possible to do the same type of calculations for large blocks of data in parallel. In this paper we therefore investigate if GPU acceleration is a measure to further speed up the established CPU-parallel BEM solver.}

\section{Introduction}
\label{author:introduction}
Every high-voltage device has to pass dielectric type tests, in which a
large voltage is applied to the device.
The test is passed if no dielectric breakdown occurs.
A breakdown usually starts from an electrode-surface with high dielectric
stress, and then propagates through the volume along a field-line of the electric field $\VE$ towards the opposite
electrode, see Fig.~\ref{Fig:Breakdown}. The propagation stops if the electric field along this breakdown
path $\gamma$ is not strong enough. For details see \cite{VHVlab}.
An inception of a streamer, i.e. the initial state of a breakdown only occurs if the criterion
\begin{equation}\label{Breakdowncrit}
  \int_{\gamma}\alpha_{\mathrm{eff}}(|\VE|) ds > K_{\mathrm{str}}
\end{equation}
is fulfilled. Here $\alpha_{\mathrm{eff}}$ is the effective ionization function that depends on the strength of the electric field $|\VE|$, and $K_{\mathrm{str}}$ is the (gas-specific) streamer constant. The prediction of a dielectric breakdown during a type test relies on the evaluation of this criterion along the most probable breakdown paths. It requires the computation of the electric field at all surface points and along field lines in the volume.

Simulation-based dielectric design became a standard procedure because a user-friendly, i.e., fast, robust, reliable, and easy-to-use computational method was developed, see \cite{Polopt1}. In the following we will first describe this boundary-element-based method and then introduce how general-purpose graphics processing units (GPGPUs) can be used to massively reduce computing times.\\
\begin{minipage}{0.49\textwidth}
\centering
\includegraphics[width=0.95\textwidth]{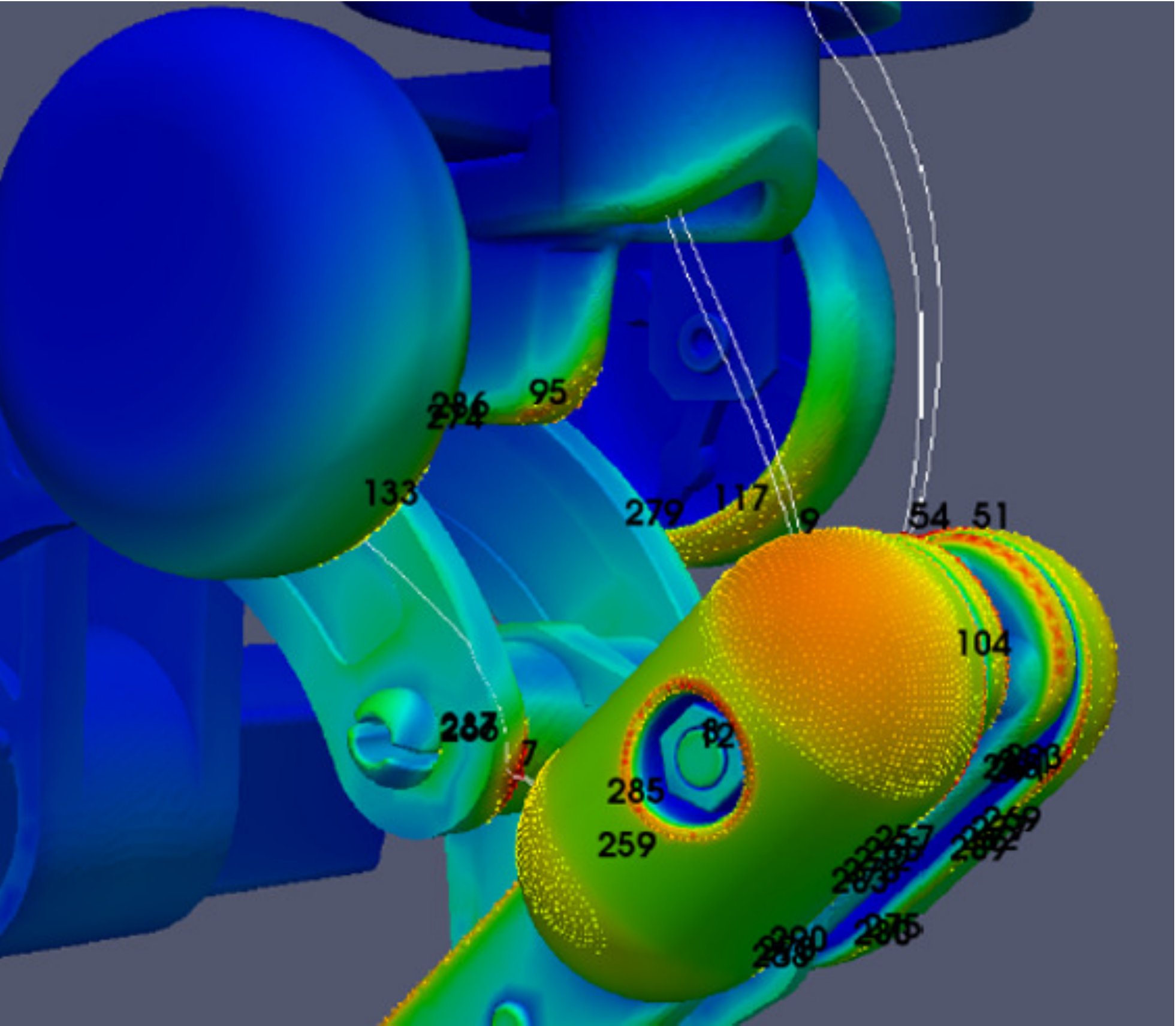}
\captionof{figure}{The electric field strength on the surface of a disconnector and possible breakdown paths along the field lines. } \label{Fig:Breakdown}
\end{minipage}
\hspace{0.02\textwidth}
\begin{minipage}{0.49\textwidth}
\includegraphics[width=0.95\textwidth]{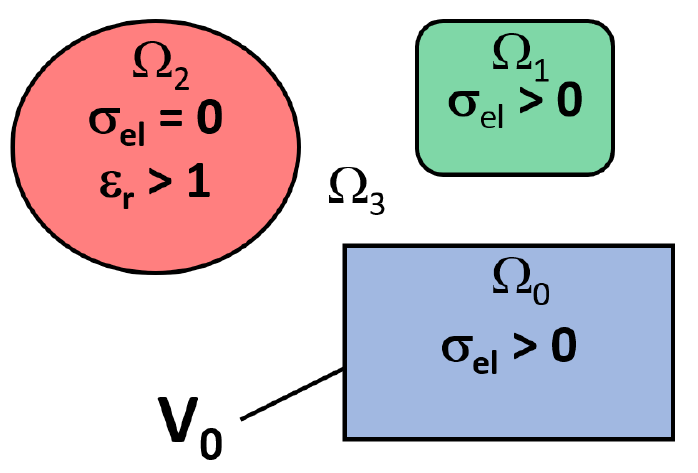}
\captionof{figure}{Example of a device that consists of a conductor $\Omega_0$ that is connected to a prescribed electric potential $V_0$, a floating conductor $\Omega_1$, an insulator $\Omega_2$ , and the unbounded exterior domain $\Omega_3$.} \label{Fig:Example}
\end{minipage}

\vspace{-0.5cm}
\section{BEM Formulation}
\vspace{-0.5cm}
In this section we derive the BEM-formulation as it is in use since many years at
ABB, see \cite{Polopt1, Polopt2}.
The device consists of the subdomains
$\Omega_0,\ldots,\Omega_{m-1}$, and the unbounded exterior subdomain is  $\Omega_{m}=\mathbb{R}^3\setminus \bigcup_{k=\{0..m-1\}} \Omega_k$, see the example in Fig.~\ref{Fig:Example}. The electric field $\VE=- \grad \varphi$ is
calculated by solving the Laplace equation
\begin{equation}\label{Laplace}
  \Div  \epsilon \grad  \varphi = 0
\end{equation}
for the electric scalar potential $\varphi$ in each of these subdomains. The permittivity is denoted by $\epsilon$.
We use an indirect formulation with a single-layer potential
\begin{equation}\label{SL}
  \varphi(\Vx) = \Psi_{SL}[\sigma](\Vx)
  = \int_{\partial\Omega} \frac{\sigma(\Vy)}{4\pi |\Vx-\Vy|} \; dS_y\,,
\end{equation}
and search for the unknown scalar virtual surface charge density $\sigma$ that is related to the physical surface charge density $\sigma_s$. Each \textbf{conductor}, i.e. each separated conducting part with electrical conductivity $\sigma_{el} > 0$, is on a constant electrical potential. If a conductor is connected to an electric potential $V_0$, like $\Omega_0$ in Fig.~\ref{Fig:Example}, then it holds
\begin{equation}
\label{conductor}
\varphi(\Vx)=V_0 \quad\forall \Vx\in \Omega_0.
\end{equation}
The electric potential $V$ of \textbf{floating conductors} like $\Omega_1$ in Fig.~\ref{Fig:Example} is unknown
\begin{equation}
\label{floatingconductor}
\varphi(\Vx)=V \quad\forall \Vx\in \Omega_1,
\end{equation}
and is to be determined by a charge neutrality condition, see \cite{Floating}. The total charge $Q$ of the floating conductor can be derived from the Gauss law as
\begin{equation}
\label{Gauss}
Q=\int_{\partial\Omega_1} \sigma_s\, dS=\int_{\partial\Omega_1} \VD \cdot \Vn\, dS.
\end{equation}
The normal component of the displacement field $\VD$ can be expressed as
\begin{equation}
\label{Jump}
\VD\cdot \Vn = \varepsilon^+ \VE\cdot \Vn = -\varepsilon^+ \grad \varphi \cdot \Vn = -\varepsilon^+ \grad \Psi_{SL} [\sigma] \cdot \Vn.
\end{equation}
Here $\varepsilon^+$ denotes the permittivity of the exterior domain.
The Neumann trace of the single layer potential and can be expressed with help of the adjoint double layer $K'$
\begin{eqnarray}
\label{ADL}
&&\grad \Psi_{SL} [\sigma] \cdot \Vn =  \frac{1}{2} \sigma + K' \sigma\mbox{ , with}\\
&&K'(\sigma)(\Vx) = \int_{\partial\Omega} \frac{ \Vx -\Vy }{4 \pi  | \Vx -\Vy |^3} \cdot \Vn(\Vx) \sigma(\Vy) dS_y\,.
\end{eqnarray}
Combining the equations \eqref{Gauss}-\eqref{ADL}  with $Q=0$ due to charge neutrality yields
\begin{equation}
\label{Eq:constraint}
\int_{ \partial \Omega_1} \frac{1}{2} \varepsilon^+ \sigma (\Vy) + \varepsilon^+ (K'\sigma)(\Vy) dS_y  = 0.
\end{equation}
We model thin floating conductive sheets only by a single surface. Then the electric fields from both sides $(\pm )$ need to be considered for charge neutrality, since
\begin{eqnarray}
\label{eq:surfcharge}
&&\sigma_s=\Vn \cdot( \VD^+ -\VD^-)\;\Longrightarrow\\
\label{Eq:sheets}
&&\int_{\partial \Omega_1} \frac{1}{2} (\varepsilon^++\varepsilon^-) \sigma (\Vy) + (\varepsilon^+-\varepsilon^-) (K'\sigma)(\Vy) dS_y  = 0.
\end{eqnarray}
There is no surface charge on \textbf{non-conductors}: $\sigma_s=0$ on $\partial\Omega_2$
\begin{equation}
\label{nonconductors}
\frac{1}{2}(\varepsilon^+ + \varepsilon^-) \sigma(\Vx) + (\varepsilon^+ - \varepsilon^-)(K'\sigma)(\Vx) = 0 \quad\forall \Vx\in\partial \Omega_2.
\end{equation}
So for our simple but quite general example of Fig.~\ref{Fig:Example} we have to solve the following set of equations:
\begin{eqnarray}\label{Final1}
 \int_{\partial\Omega} \frac{\sigma(\Vy)}{4\pi |\Vx-\Vy|} \; dS_y &=& V_0 \quad\forall \Vx\in \partial\Omega_0\\
\label{Final2}
 \int_{\partial\Omega} \frac{\sigma(\Vy)}{4\pi |\Vx-\Vy|} \; dS_y-V &=& 0 \quad\forall \Vx\in \partial\Omega_1\\
\label{Final4}
\int_{\partial \Omega_1} \frac{1}{2} \varepsilon^+ \sigma(\Vy) + \varepsilon^+  (K'\sigma)(\Vy) dS_y  &=& 0\\
\label{Final3}
 \frac{1}{2}(\varepsilon^+ + \varepsilon^-) \sigma(\Vx) +  (\varepsilon^+ - \varepsilon^-)(K'\sigma)(\Vx) &=& 0 \quad\forall \Vx\in \partial \Omega_2
\end{eqnarray}
The solution of the system of equations \eqref{Final1} - \eqref{Final4} yields  the virtual surface charge distribution from which the electric field can be compute at any point in space as
\begin{equation}
\VE(\Vx) = \int_{\partial\Omega} \frac{ \Vx -\Vy }{4 \pi | \Vx -\Vy |^3} \sigma(\Vy) dS_y \quad \forall\Vx \in \mathbb{R}^3.
\end{equation}

\vspace{-1cm}
\section{Discretization}
\label{sec:Discretization}
\vspace{-0.5cm}
We use a collocation boundary element approach: the surface $\partial\Omega$ is
represented by a collection of triangles $\tau_1,\ldots,\tau_N$
with vertices $\Vx_1,\ldots,\Vx_n$.
The unknown function $\sigma$
\begin{equation*}
  \sigma_h(\Vy) = \sum_{j=1}^n u_j \psi_j(\Vy),
\end{equation*}
is approximated by suitable basis functions $\psi_1,\ldots,\psi_n$. This approach is to be inserted into \eqref{Final1}-\eqref{Final4} and this set of equations  is only required to hold in the
\emph{collocation points} $\Vx_1,\ldots,\Vx_n$.
This is an $n+N_{fl}$-dimensional system of linear equations, with $N_{fl}$ being the number of floating conductors. The implementation of this approach poses a number of challenges:
\begin{itemize}
\item High-voltage devices have smooth curved surfaces in order to avoid field enhancements. We use piecewise quadratic parametrizations with (curved)
  triangles $\tau_1,\ldots,\tau_N$ to minimize the geometrical discretization error.
\item The entries of the matrix corresponding to the linear
  system have to be computed.
  In the established CPU-based method we employ an MPI-parallel implementation of suitable
  quadrature rules.
\item The system of linear equations has to be solved.
  We use Krylov subspace methods, since these methods only
  need matrix-vector multiplications, which can easily
  provided in the established method by the MPI-based distributed representation of
  the matrix.
\end{itemize}

\vspace{-1cm}
\section{GPGPU Quadrature}
\label{sec:GPGPU}
\begin{figure}[htb]
\begin{center}
\includegraphics[width=0.95\textwidth]{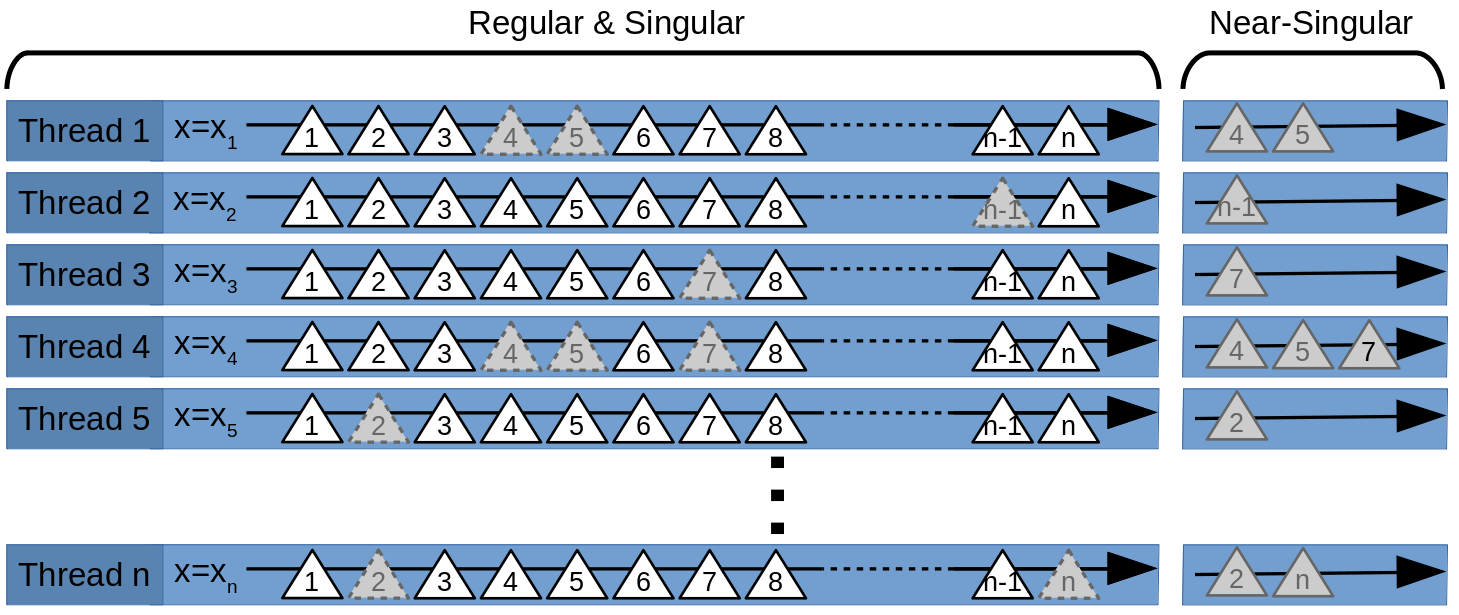}
\caption{Near-singular assembly: The gray triangles are near-singular and are computed after the regular and singular triangles. } \label{Fig:near}
\end{center}
\end{figure}

\vspace{-0.5cm}
Even on modern processors of parallel computers it is time-consuming to compute the matrix entries $v_{ij}$  for complex industrial geometries. It is advantageous to calculate the surface integral for pairs of collocation points and triangles.
\begin{align}\label{Matrix}
  v_{ij} &= \int_{\partial\Omega} \frac{\psi_j(\Vy)}{4\pi |\Vx_i-\Vy|} \;dS_y &
  &\text{ for all } i,j\in\{1,\ldots,n\}
\end{align}
Matrix assembly is however ideally suited for parallelization on SIMD-type of processors, because
the integrations require us to
carry out mostly identical operations for \emph{all} matrix entries.
The currently most common processors are general-purpose
graphics processing units (GPGPUs) like the NVidia Tesla\texttrademark{} or
AMD Radeon Instinct\texttrademark{} cards that contain thousands of
floating-point arithmetic units and offer teraflop-level performance.
Porting the quadratures to GPGPUs poses challenges:
\begin{itemize}
  \item The most powerful GPGPU models are equipped with fast local memory. We have to ensure that geometrical data is efficiently transferred to the local memory.
  \item We are using piecewise (mapped) linear basis functions and multiple
triangles contribute to the same matrix entry. Thus we have to avoid
collisions between different triangles that may try to simultaneously update the same matrix entry.
\end{itemize}
For the parallelization it is beneficial that each collocation point corresponds to one row in the matrix. In our implementation we take advantage of this fact by assigning each collocation point to a thread. Then we consecutively iterate through all triangles of the discretization and simultaneously compute the  contribution of a triangle to all collocation points, see Fig.~\ref{Fig:near}. This way collisions can be avoided. For good performance it is required that all threads execute exactly the same operations. This can be ensured by grouping the collocation points depending on whether they are on the surface of a conductor or part of a dielectric interface.

For the quadrature we distinguish between regular, near-singular, and singular integration of the pairs of triangles and collocation points. If the collocation point is a node of the triangle then it is a singular pair. When the distance $D$ between the circumcenter of the triangle and the collocation point is larger than a threshold that scales with the circumradius $R$ of the triangle like
\begin{equation}
D > \eta\cdot R
\end{equation}
then it is a regular pair. We used the scaling parameter $\eta = 1.2$. All other pairs are near-singular. Circumcenters and circumradii of the triangles can be precomputed.

All three types of pairs can be integrated by using the well-known Duffy-transformation, see \cite{DUFFY}.
 This is straightforward for regular and singular pairs, only the near-singular pairs need to be integrated adaptively for accuracy. They frequently occur, e.g. in cases with narrow gaps in the geometry, or during the postprocessing, when a point near the surface is to be evaluated. In this near-singular case we first compute the point of the triangle that is closest to the collocation point. Next we subdivide the triangle into smaller triangles
  such that this closest point is a corner-node of a subdividing triangle.
  Then we  again employ the Duffy-transformation to integrate over all smaller triangles.
  This yields an adaptive quadrature with increased accuracy around the closest point. The adaptivity can impact the performance of a GPU-computation badly if no attention is paid, because then there is divergence in the control flow on the GPU. In order to minimize this divergence, we first compute the regular and singular pairs, and deal with the near-singular integrals later.

The categorization into regular, near-singular and singular pairs is carried out on the fly during the iteration through the triangles.
 If a pair is marked as near-singular, then it will be marked as not processed, see Fig.~\ref{Fig:near}. They are computed in parallel by using the subdivision method after the regular and singular cases have been completed.  This strategy allows that all three types of integrations are carried out in parallel, without the need to mix operations.
 \vspace{-5mm}
 \vspace{-7mm}

 The full matrix may not fit in the memory of one GPU for larger problems. Therefore, but also to speed up the computation we use multiple GPUs. Due to the independence of matrix rows we split the matrix into multiple blocks of rows that can be computed and stored independently on different GPUs.

%
%
%

\vspace{-0.5cm}
\section{Numerical Experiments}
\vspace{-0.5cm}
In this chapter we show some examples that were computed with the novel GPU-implementation that is based on the H2Lib package, see \cite{H2Lib}. We first validated our implementation for an axial-symmetric case. We compared the results of the H2Lib with the results of the already existing simulation tools Polopt (3D) see \cite{Polopt2},  and Elfi (2D) see  \cite{Polopt1}. Next we compared the performance of the new GPU-parallel H2Lib implementation with the performance of the existing MPI-parallel Polopt tool.

\vspace{-0.5cm}
\subsection{Validation}
\vspace{-0.5cm}
The benchmark problem that is used to validate the GPU-implementation in H2Lib is a bushing, see Fig.~\ref{Fig:bushing}. It consists of an insulator that is wrapped around a conductor on high-voltage. Five thin conducting sheets are embedded in this insulator. They accomplish the field grading. The outermost is grounded, and the potentials of the other (floating) sheets are unknown and are to be determined. The sheets are treated as single surfaces according to equation \eqref{Eq:sheets}. Their potentials were computed with all three solvers. POLOPT and H2Lib use the same mesh with 3'526 nodes. The results agree very well, see Table~\ref{tab:bushing}. The small remaining differences are due to the use of different quadratures.

\begin{minipage}{0.5\textwidth}
%
\begin{tabular}{|l|l|l|l|}
\hline
 & ELFI(2D) & POLOPT & H2LIB\\\svhline
V\textsubscript{foil1} & 70.8kV & 70.15kV & 70.22kV \\
V\textsubscript{foil2} & 51.4kV & 50.47kV & 50.50kV \\
V\textsubscript{foil3} & 35.0kV & 34.00kV & 34.02kV \\
V\textsubscript{foil4} & 18.9kV & 18.02kV & 18.02kV \\ \hline
\end{tabular}
\captionof{table}{Potentials of the conducting sheets of the bushing}
\label{tab:bushing}
\end{minipage}
\hspace{1mm}
\begin{minipage}{0.5\textwidth}
\centering
\includegraphics[width=0.95\textwidth]{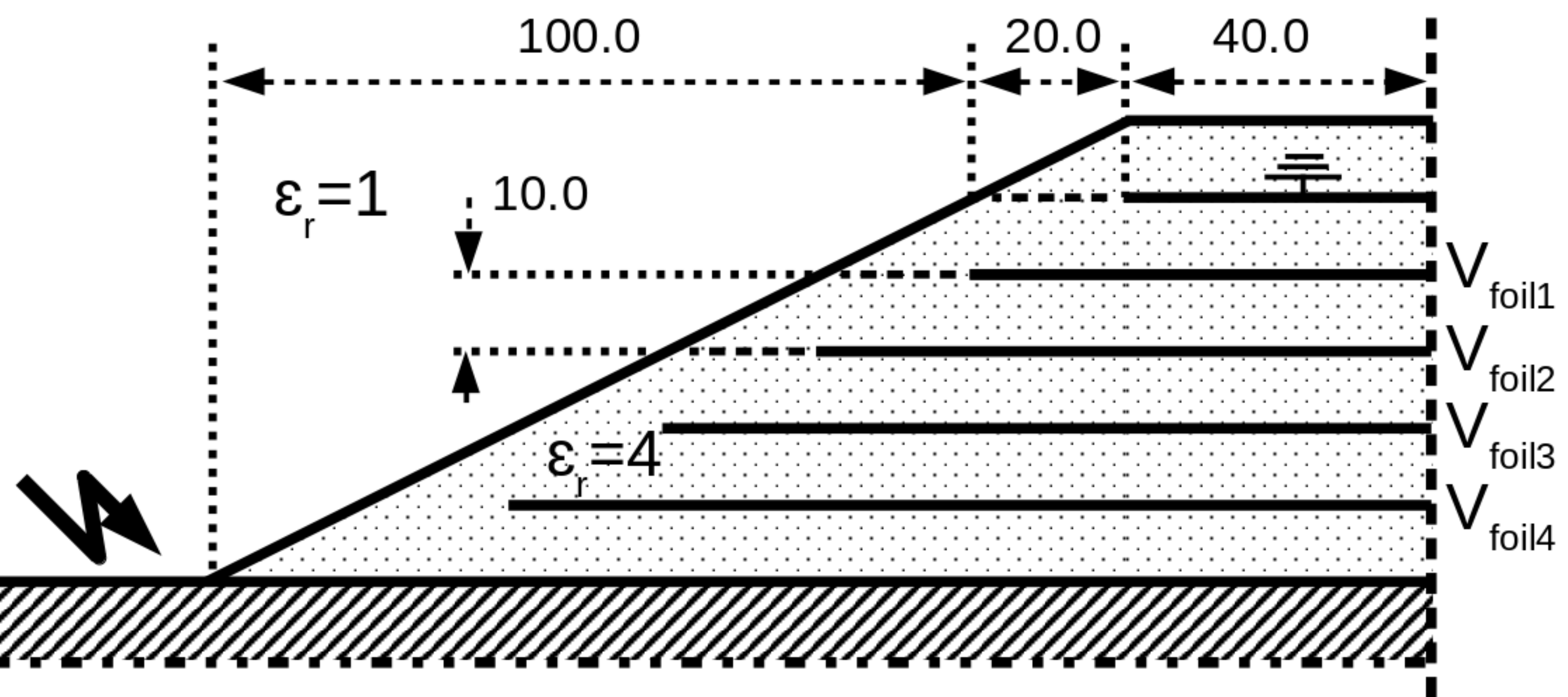}
\captionof{figure}{Cross-section of the bushing}
\label{Fig:bushing}
\end{minipage}

\vspace{-0.5cm}
\subsection{GPU-Acceleration}
\vspace{-0.5cm}
After the successful validation of the H2Lib implementation, we compared the computing times for the GPU-parallel H2Lib and the CPU-parallel Polopt. In both cases we assembled the dense BEM matrix and solved the system with an iterative GMRES with diagonal preconditioner. So the expected numerical work  is quadratically depending on the degrees of freedom, i.e. here the number of nodes (\# Nodes).

\begin{minipage}{0.49\textwidth}
\hspace{-2mm}
\includegraphics[width=1.05\textwidth]{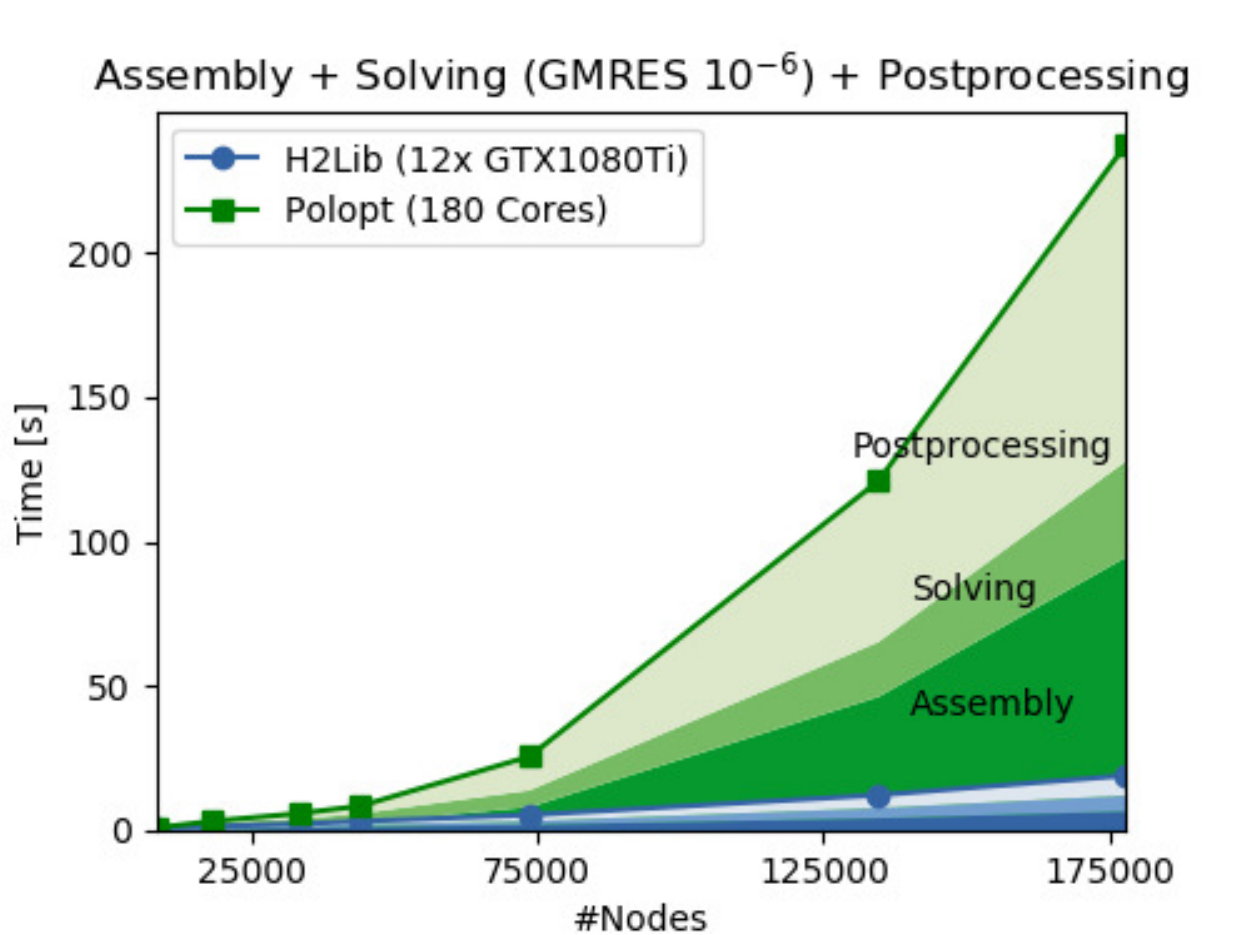}
\captionof{figure}{Cumulative times for assembly, solving and  surface electric field computation for POLOPT and H2LIB.}
\label{Fig:Polopttime}
\end{minipage}
\hspace{0.01\textwidth}
\begin{minipage}{0.49\textwidth}
\hspace{-2mm}
\includegraphics[width=1.05\textwidth]{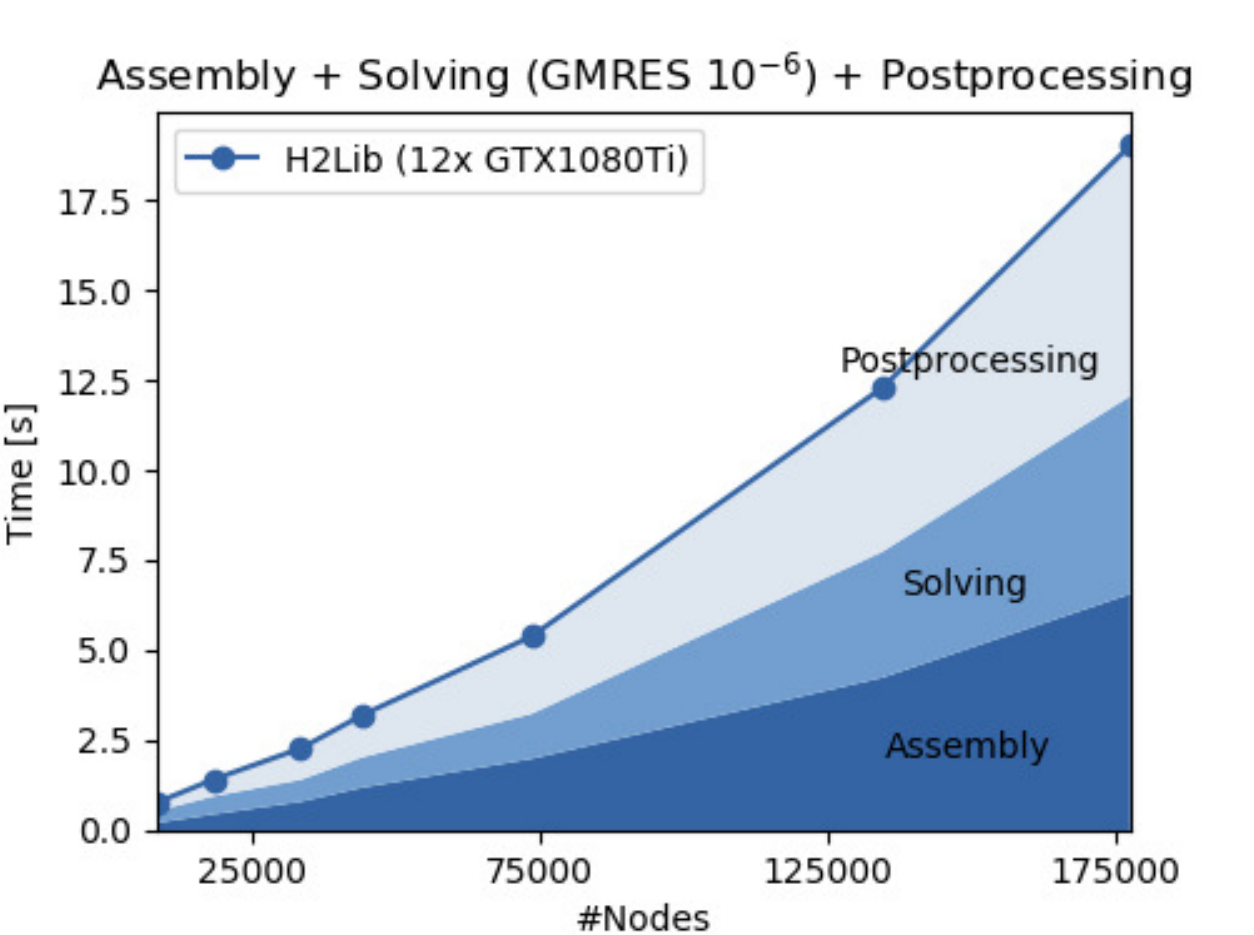}
\captionof{figure}{Cumulative times for assembly, solving and  surface electric field computation for the H2Lib only.}
\label{Fig:times}
\end{minipage}
\vspace{2mm}

Figure \ref{Fig:Polopttime} shows the times that it took for matrix-assembly, iterative solution, and computation of the electric field at the surfaces of a realistic high-voltage device for different mesh-sizes. POLOPT used  180 CPU cores distributed over 5 nodes with two 18-core CPUs each. The CPU cluster was optimized for these calculations because it is in use by ABB product designers. For H2Lib we used a total of 12 NVidia GTX 1080Ti, distributed over multiple nodes. For POLOPT we clearly recognize the quadratic scaling.
The H2Lib also scales quadratically, see Fig.~\ref{Fig:times}, the proportionality constant seems however to be much better than the one for Polopt.\\

\noindent
 \begin{minipage}{0.49\textwidth}
\begin{tabular}{ | r | r | r | r |  }
\hline
 \,\#Nodes\, & \,\#GPUs\, & \,H2Lib\, & \,POLOPT\, \\ \svhline
 68'218 & 2		& 18s	& 33s \\
140'183 & 8     & 19s   & 130s \\
232'029 & 20    & 26s   & 371s \\
330'706 & 40    & 34s  	& 702s \\
432'084 & 72	& 41s   & 732s$^*$ \\ \hline\noalign{\smallskip}
\end{tabular}
\\
$^*$ 360 Cores
\vspace{-1mm}
\captionof{table}{Calculation times for different meshes. H2Lib used multiple GTX1080Ti, POLOPT is executed on 180 CPU-Cores.}
\label{Tab:exk0gis}
\end{minipage}
\hspace{0.02\textwidth}
\begin{minipage}{0.49 \textwidth}
\begin{tabular}{ | l |  p{1.8cm} |  p{2cm} | }
\hline
 & POLOPT    &  H2Lib  \\
 & 180 Cores &  4xTesla P100\\
\svhline
 Assembly & \multicolumn{1}{r|}{11s \hspace{3mm}  } & \multicolumn{1}{r|}{2s \hspace{1mm} }\\
Solving & \multicolumn{1}{r|}{11s \hspace{3mm}  } & \multicolumn{1}{r|}{2s \hspace{1mm} }\\
$\VE$ surface & \multicolumn{1}{r|}{11s \hspace{3mm} } & \multicolumn{1}{r|}{2s \hspace{1mm} } \\ \hline
Total & \multicolumn{1}{r|}{33s \hspace{3mm}  } & \multicolumn{1}{r|}{6s \hspace{1mm} }\\\svhline
Evaluation & 102s for  & 67s for \, 9219\\
            & 17 fieldlines &   fieldlines \\ \hline
\end{tabular}
\vspace{-3mm}
\captionof{table}{Time for computation and fieldline evaluation. Model size: 68'218 Nodes. Polopt in serial for the evaluation.}
\label{Tab:fieldlines}
\end{minipage}\\

\noindent
We computed another larger example see Table \ref{Tab:exk0gis}.  We used 180 cores for POLOPT for all meshes except the largest one, where we used 360 cores. The number of GPUs were chosen such that the matrix  fitted in the combined memory of all GPUs in single precision for H2Lib. Again the GPU-parallel implementation clearly outperforms the CPU-parallel version. We also evaluated the breakdown criterion along the most probable breakdown paths (field lines), see Table \ref{Tab:fieldlines}. The acceleration seems even higher, however the evaluation is only implemented in serial in Polopt.

%

\vspace{-0.5cm}
\section{Conclusions and Outlook}
\vspace{-0.3cm}
The usage of GPUs drastically accelerates the computation of electrostatic fields, as well as their evaluation with respect to breakdown inception.
This opens the door not only for the computation of larger problems, but also for the inclusion of additional physical models (e.g. for surface charging), or for optimization.
Most promising is the combination of the GPU-acceleration with a  compression technique, see \cite{GCA,GCA_SCEE}. Another strong acceleration is to be expected in this case. Moreover the asymptotical behavior will no longer depend quadratically  on the degrees of freedom.



\end{document}